\documentclass{article}
\usepackage{graphicx} 

\usepackage{xcolor} 
\usepackage[a4paper, total={6in, 8in}]{geometry}
\usepackage{amssymb,amsfonts,amsthm,amsmath}
\usepackage[toc]{appendix}
\usepackage{hyperref} 
\usepackage[sanserif,full]{complexity}
\usepackage{mathtools}
\usepackage[color=blue!25]{todonotes}
\usepackage{float}
\usepackage{enumitem}
\usepackage{subcaption}

\title{A Mathematical Murmuration}
\author{Mathison Knight\footnote{Center of Mathematics and Theoretical Physics, University of Hertfordshire, UK.} \,\ and  Shanti Pise\footnote{School of Physics, Engineering and Computer Science, University of Hertfordshire, UK.} }
\date{}

\begin{document}

\maketitle

\textit{We write this article as a summary of a large group movement activity as part of a STEAM (Science, Technology, Engineering, Arts and Mathematics) event at the University of Hertfordshire.}\\

\noindent Intuition of the body, and intuition of mathematics feel like two sides of the same coin, but are often ignored, so how can we bring these together? Dance is a natural, human activity, which has been observed in almost all human societies, as discussed in \cite{Hannay2024Origins}. The embodiment of mathematics in dance (or more generally, movement and motion) and vice versa has been studied in a wide range of environments, e.g., \cite{Nasiakou2019Maths,watson2005dance,96943352-68c9-333f-bb79-ed8762d47f90,gerofsky2010mathematical}. The main focus in these is on the individual and their experience, or at most in small groups (\cite{vogelstein2019reenacting}) but we are social beings, and both dance and mathematics are often not done in near isolation... It is then natural to ask how can the embodiment of mathematics be brought to large groups?

We explore this by bringing together 377 Year 4 Students (8/9 year olds), albeit in two separate groups, to create a large group mathematical dance routine. Our aim was to demonstrate embodiment, not just with the students by themselves, but as a joint collaborative effort- as much mathematics is. In such we hoped to increase the students awareness of shape and space, both mathematically and physically, by having them perform timed movements to music.

Our main goals and questions were the following:

\begin{itemize}
    \item We wish to look at how large group choreography can help to develop mathematical reasoning.
    \item We wish to see how students (and ambassadors) experience large group activities, with the aim of showing the participants that every individual matters in a routine.
    \item We wish to show mathematics as a creative pursuit that is not necessarily procedural or rule based.
\end{itemize}
\large \textbf{ Can we assist students in viewing mathematics as a creative and collaborative pursuit with large group choreography?}
\normalsize

\section{In the Beginning...}
We invite the reader to join us on the journey of the session we conducted.

After a brief introduction of ourselves, we showed the audience a video of starlings, making beautiful, fluid patterns that are highly coordinated, and seemingly artificial. \cite{Claireonline} 

The students were shocked at how this is done, and we guided them into asking questions; What is this? What are the birds doing? Are they communicating with each other? Is this real or made by artificial intelligence? Are they doing mathematics, or dance, or both?

The concept is known as a murmuration, and is often seen in nature, we chose this as a solid foundation of bringing together previous knowledge and new information (\'a la constructivism), particularly as we were unsure of how much the students would have been familiar with, beyond the guidance from the national curriculum. A central thing we all have in common is indeed nature, and the experience of the world around us. 

Not many students had heard of this concept before, some thinking that murmuration may be linked to the other definition  of murmuring- `To complain in low muttered tones; to give voice to an inarticulate discontent; to grumble.' \cite{murmur}.

After some discussion, we found that participants were intrigued by their ability to avoid collisions whilst changing directions almost instantaneously. We opened the floor to ask why birds did this, and many thought of very intuitive answers- to evade predators, to search for food.

To solidify the idea that this is indeed a group endeavour, we asked the students whether they can make a murmuration by themselves? After some (short lived) attempts from the pupils, we landed on a firm `no'. Hence, we conjure the concept of \textbf{large group embodiment}. Embodiment shows us that the body is all you need, but how much more could we accomplish with the consideration of others in our space?

\section{Demonstration and Choreography}

Before we gave the choreography to the pupils, we introduced a group of 16 Ambassadors from the University, these 
are undergraduates, postgraduates and alumni, some of which work within the Widening Access team and others of which work in the Physics, Astronomy and Mathematics outreach programme. These staff members have backgrounds in many different fields from Law to Education to Mathematics. We wished to give the students a goal to work toward, their own murmuration which would look something similar to the Ambassadors one. 

A key observation we had was that the students were trying to figure out the movements during the demonstration, 
 from personal observation, it is clear that curious students often try to work things out before teachers have shown the solution in mathematics and this gave us a concrete visual understanding of what is happening in the mind of the pupils; this may help other educators in figuring out what pupils are thinking, when they are resistant to sharing.

 Students also asked us and ambassadors how long they had been practising for, and were visibly shocked when they were told that they had only had an hour to learn the routine! We hope that this aided them in the confidence they brought into the dance.

Below, we include the video of our ambassadors performing to the students for your viewing. 

https://youtu.be/8NWNw9hTirk
 
 We followed the demonstration with a round of reflective questioning, and  we invite you to answer the following:  What shapes can be observed? What mathematical concepts can you identify? How are patterns being formed on the individual level? On the group level?

We went through some basic movements and then handed the floor over to the pupils...

The routine that we shall describe here is what the pupils performed, first we note that this dance only
consists of moving hands and arms, and so this can be done sitting. We will run through the single movements and then summarise how these came together as a group

Beginning with closed fists, as the music springs to life, we extend our arms, pointing forward on the beat, mapping out a literal geometrical point.  We grow, in dance,  sweeping the y-axis, moving one hand up and down, and repeating with the other hand. This demonstrates expanding our points into lines, and increasing the dimension we are working in. We transition into two dimensional space with a gentle wave movement, before we look at the major shape in the routine. 

The Square.

We form this in two different ways, swiping our hands through space to give the impression of two dimensional space through one dimensional lines. 

The first way is as follows:
\begin{enumerate}[start=0]
    \item Begin with both arms flat upon each other, with the left under the right.
    \item Bring the right arm up, so the elbow is positioned next to the left hand, and it is perpendicular to the left arm. 
    \item Sweep the left arm upward, staying perpendicular to the right arm, and stop when both sets of fingers are touching.
    \item Rotate both arms by 90 degrees.
    \item Sweep the right arm downward, staying perpendicular as before. Close the left arm back down to the right arm, and repeat!
    
\end{enumerate}

This is shown in Figure \ref{sq1} below.

\begin{figure}[H]
  \centering
  \begin{minipage}[b]{0.23\textwidth}
    \includegraphics[width=\textwidth]{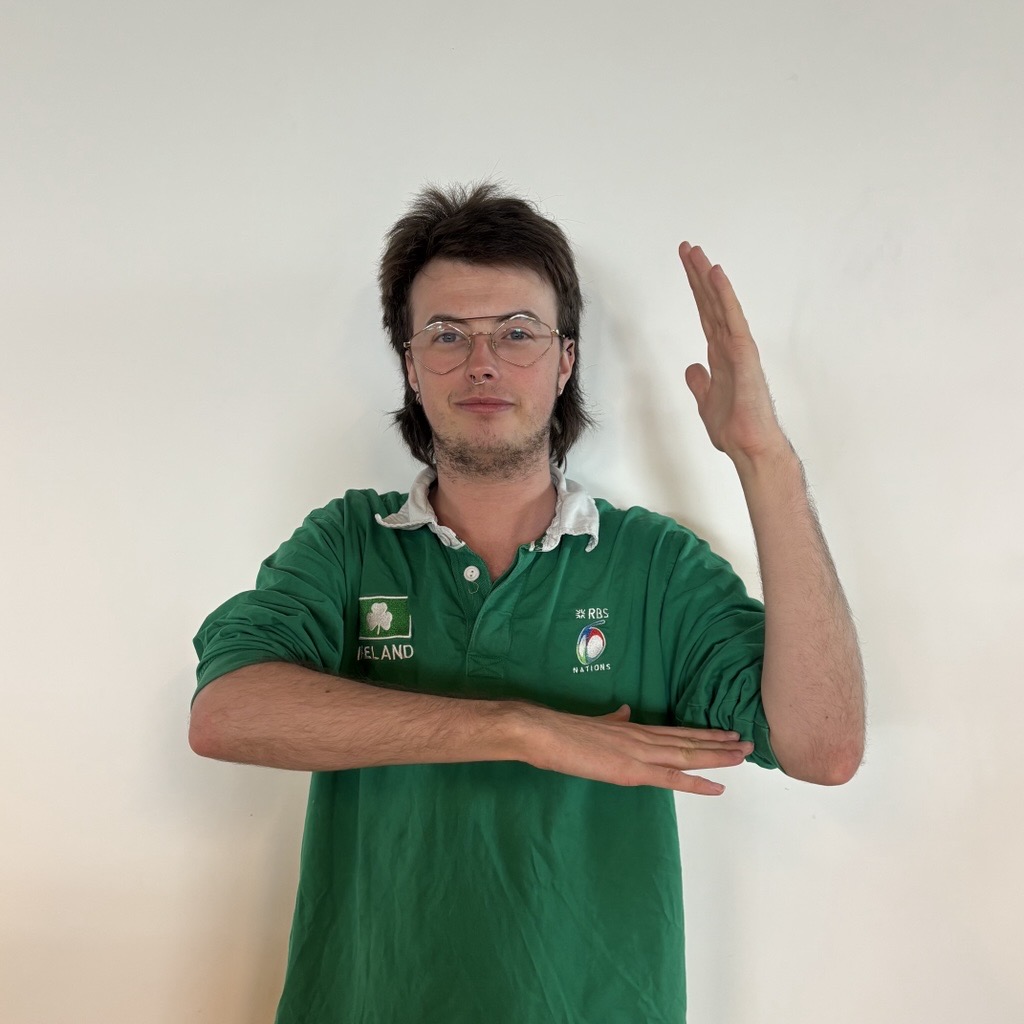}
    \subcaption{Step 1.}
  \end{minipage}
  \begin{minipage}[b]{0.23\textwidth}
    \includegraphics[width=\textwidth]{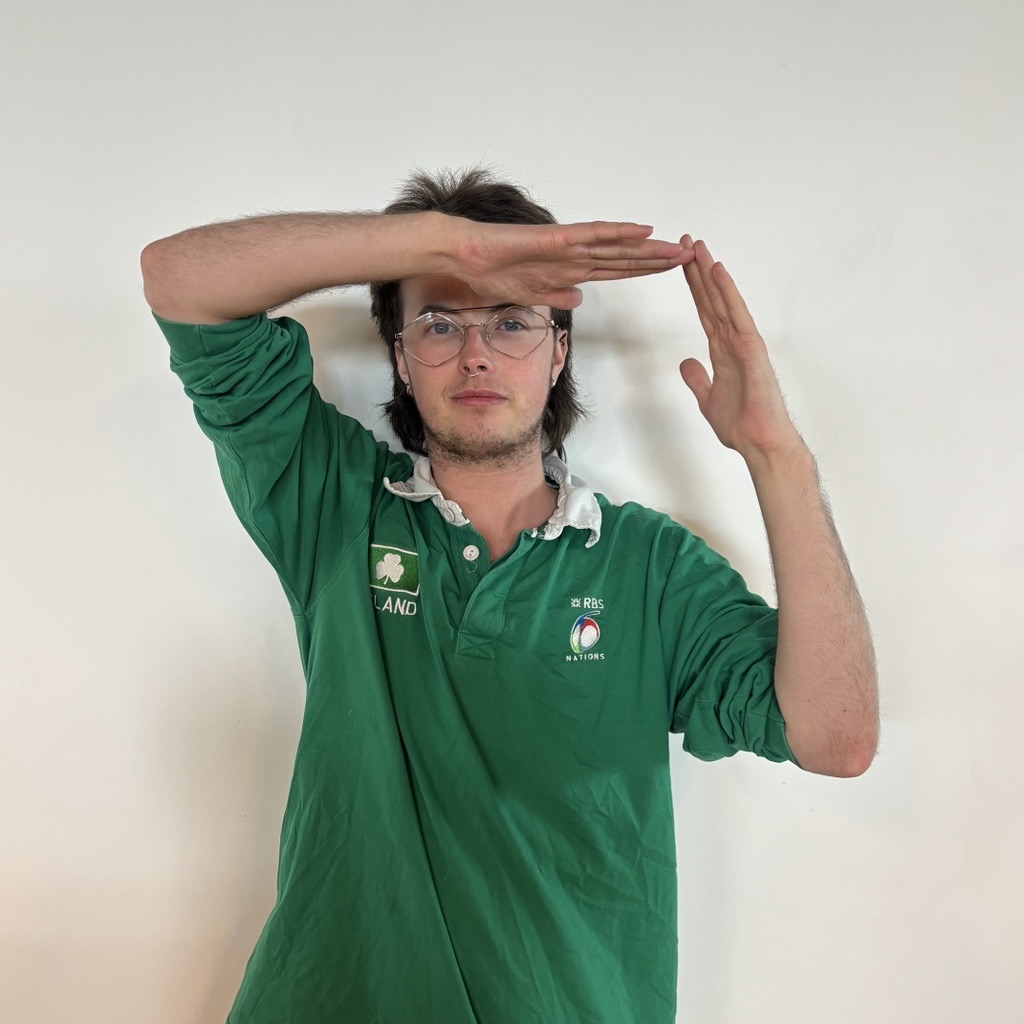}
  \subcaption{Step 2.}
  \end{minipage}
  \begin{minipage}[b]{0.23\textwidth}
    \includegraphics[width=\textwidth]{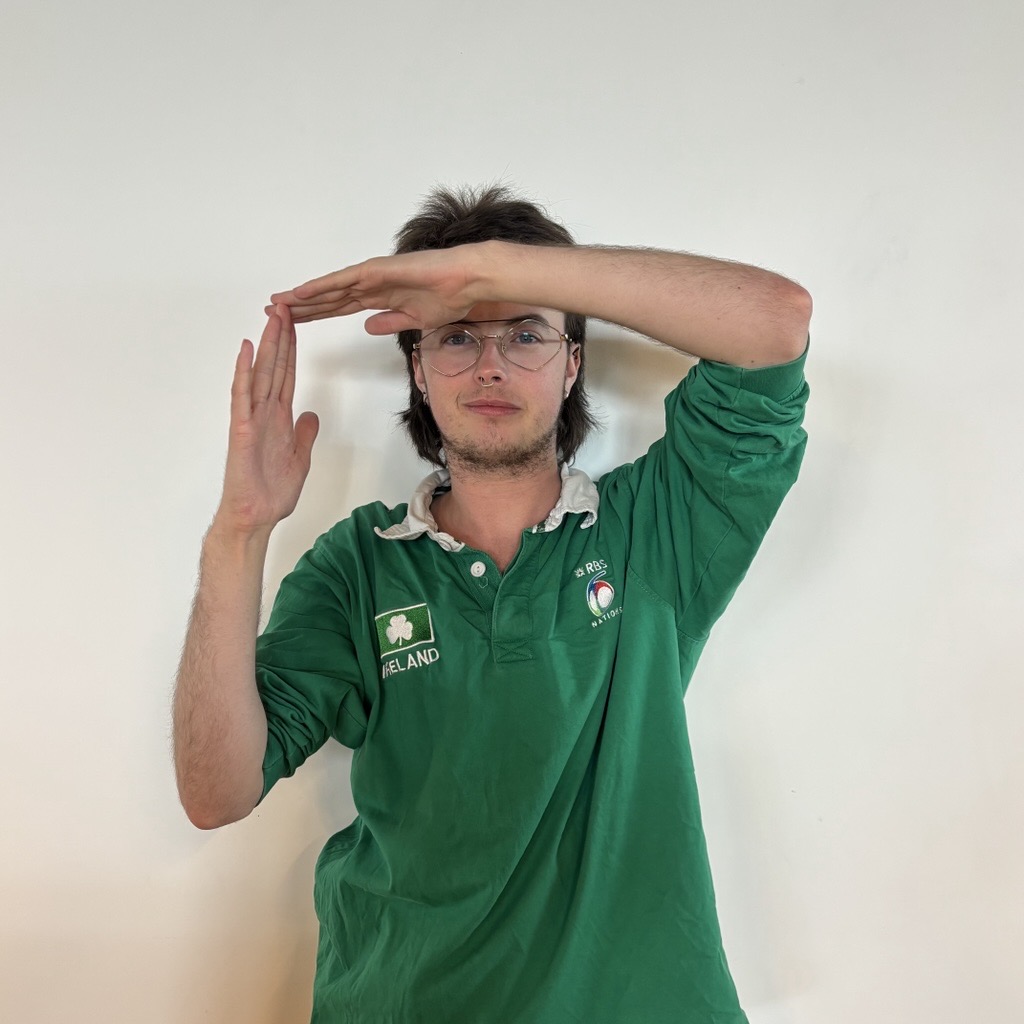}
  \subcaption{Step 3.}
  \end{minipage}
   \begin{minipage}[b]{0.23\textwidth}
    \includegraphics[width=\textwidth]{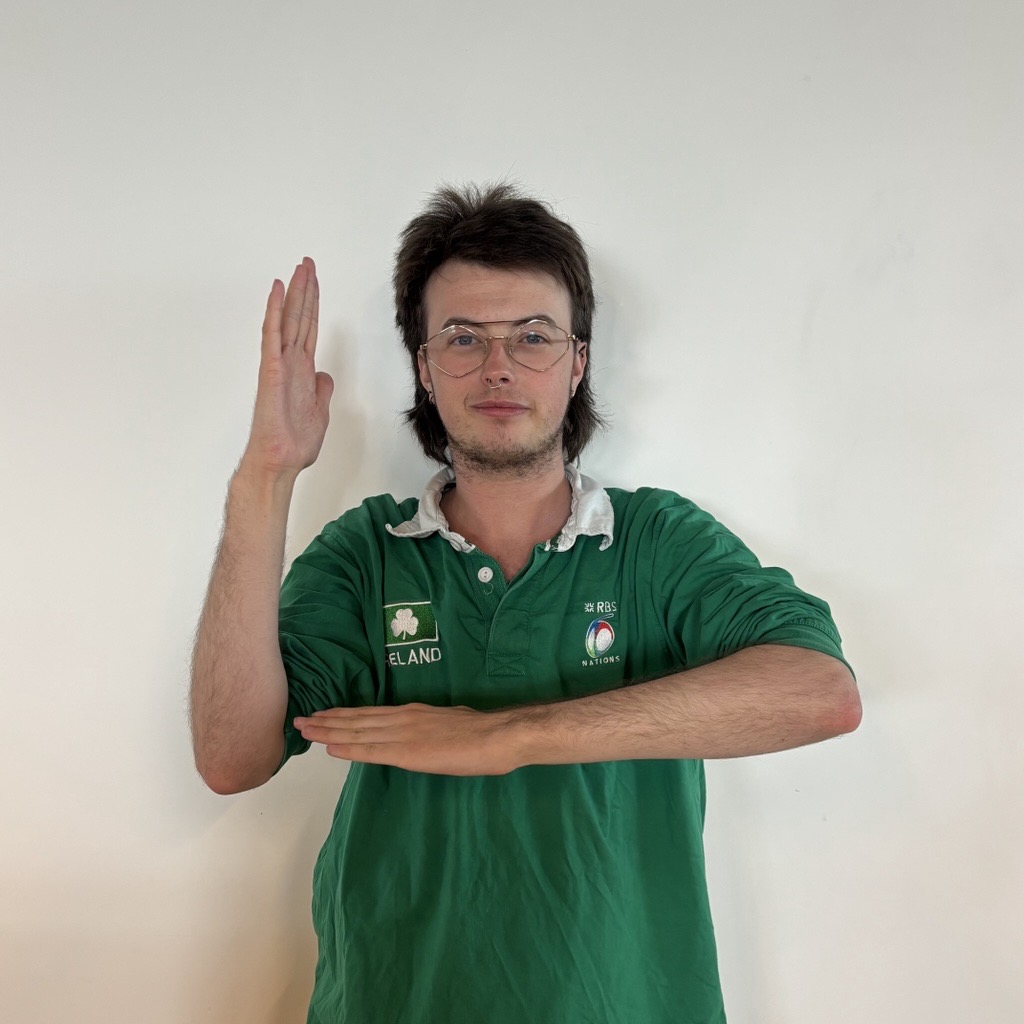}
 \subcaption{Step 4.}
  \end{minipage}
  \caption{Square Routine One}
  \label{sq1}
\end{figure}

The second iteration followed like so:
\begin{enumerate}
    \item Begin with arms crossed over the torso, the right over the left, so this looks like an X.
    \item Swipe the right arm along the left arm in a fluid motion, stopping when reaching the fingertips, so this lines are perpendicular. 
    \item Then move the left arm upwards along the right arm to make a point at forehead level, finger tips touching, showing a right angle in the square.
    \item Swipe the right arm back down along the left, until you reach about half way, and then move the left arm in a fluid motion behind the right, to reach Step 1 again!
\end{enumerate}

\begin{figure}[H]
  \centering
  \begin{minipage}[b]{0.23\textwidth}
    \includegraphics[width=\textwidth]{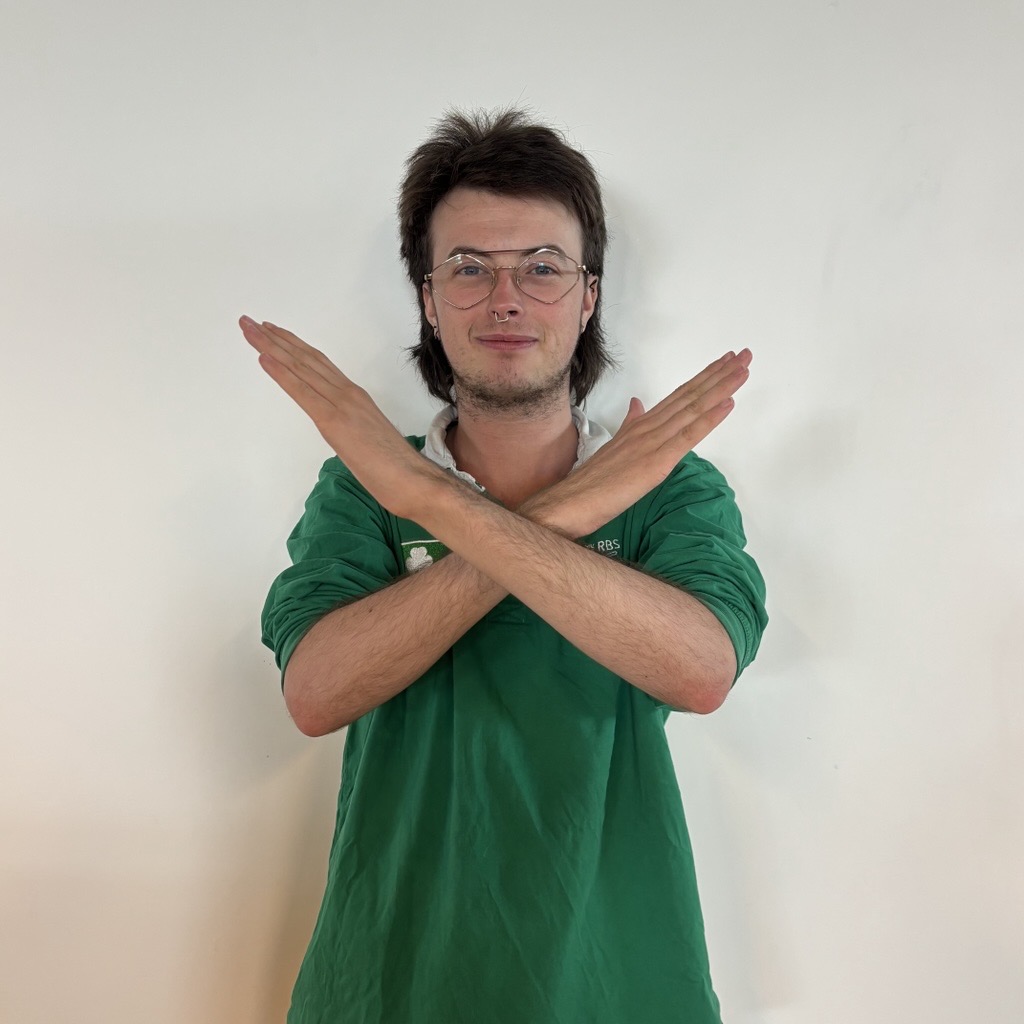}
    \subcaption{Step 1.}
  \end{minipage}
  \begin{minipage}[b]{0.23\textwidth}
    \includegraphics[width=\textwidth]{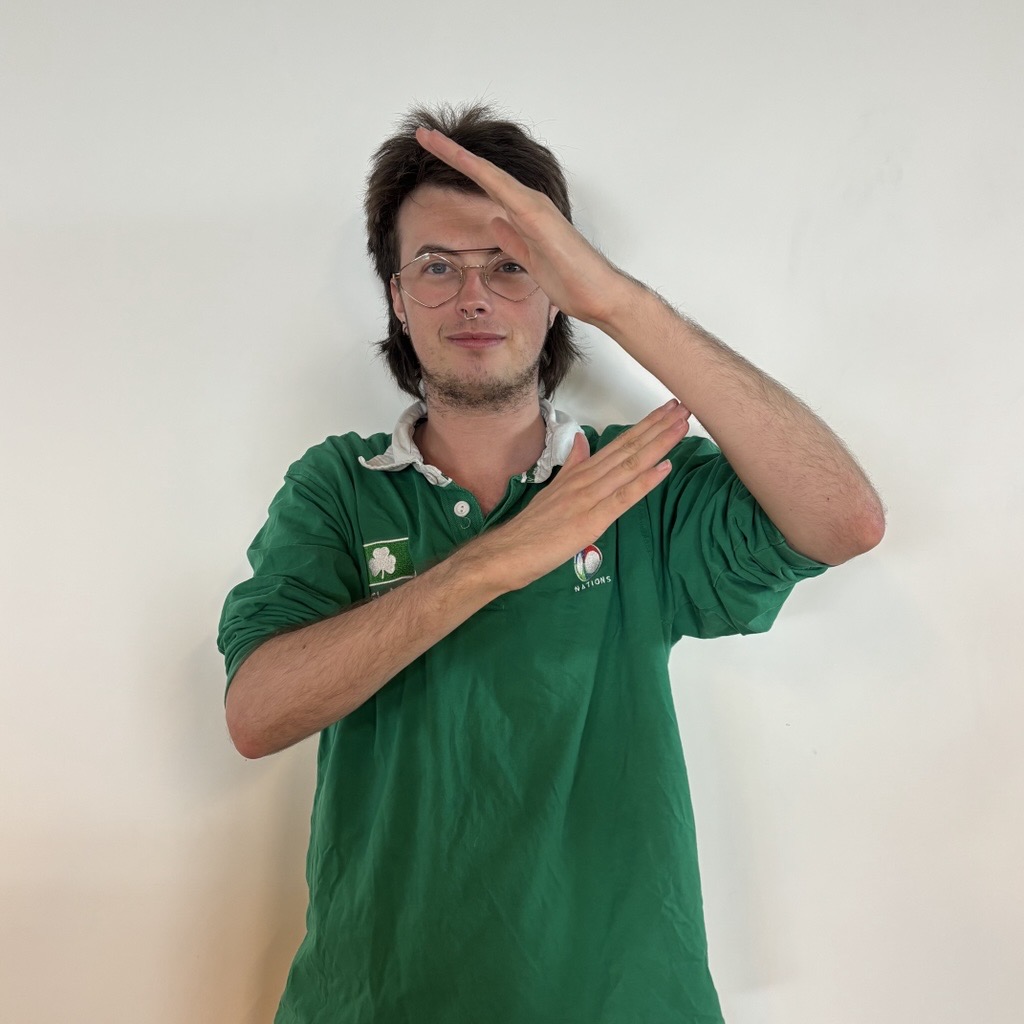}
  \subcaption{Step 2.}
  \end{minipage}
  \begin{minipage}[b]{0.23\textwidth}
    \includegraphics[width=\textwidth]{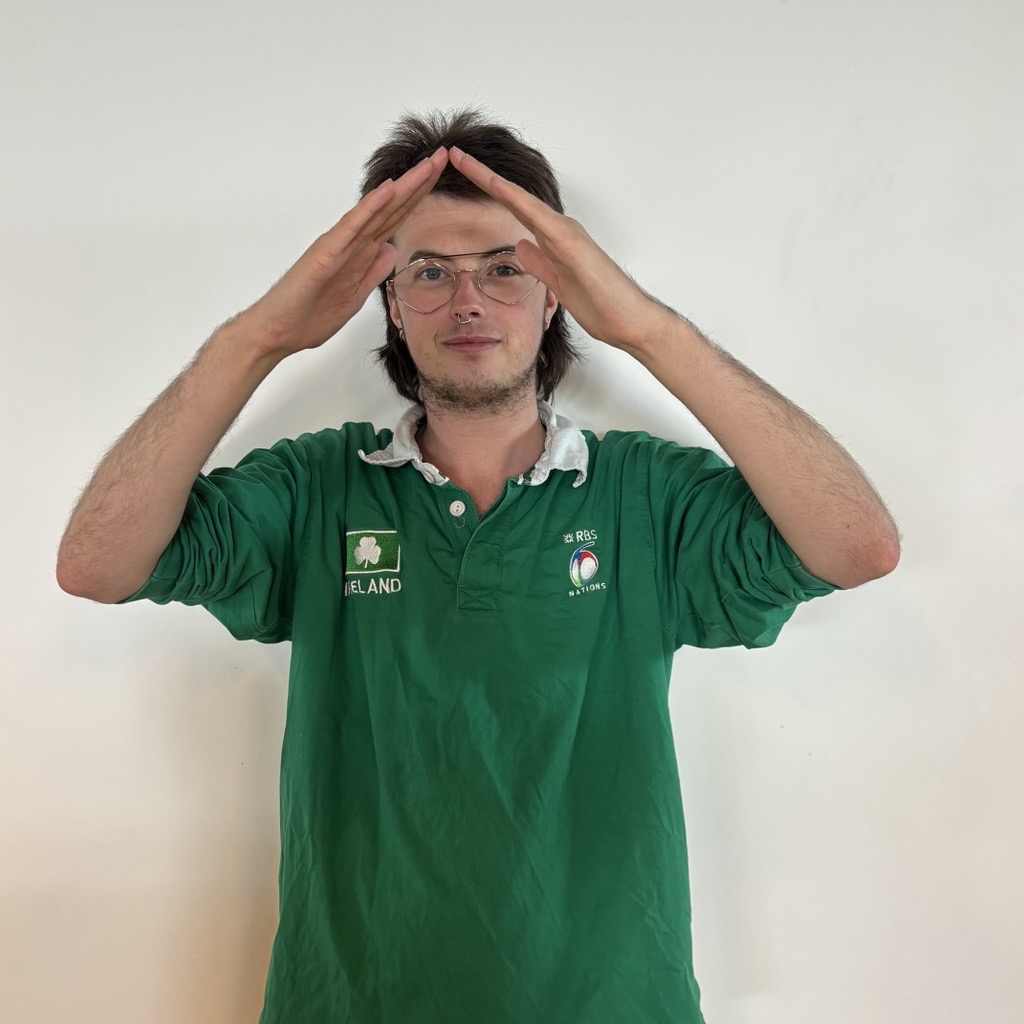}
\subcaption{Step 3.}
  \end{minipage}
   \begin{minipage}[b]{0.23\textwidth}
    \includegraphics[width=\textwidth]{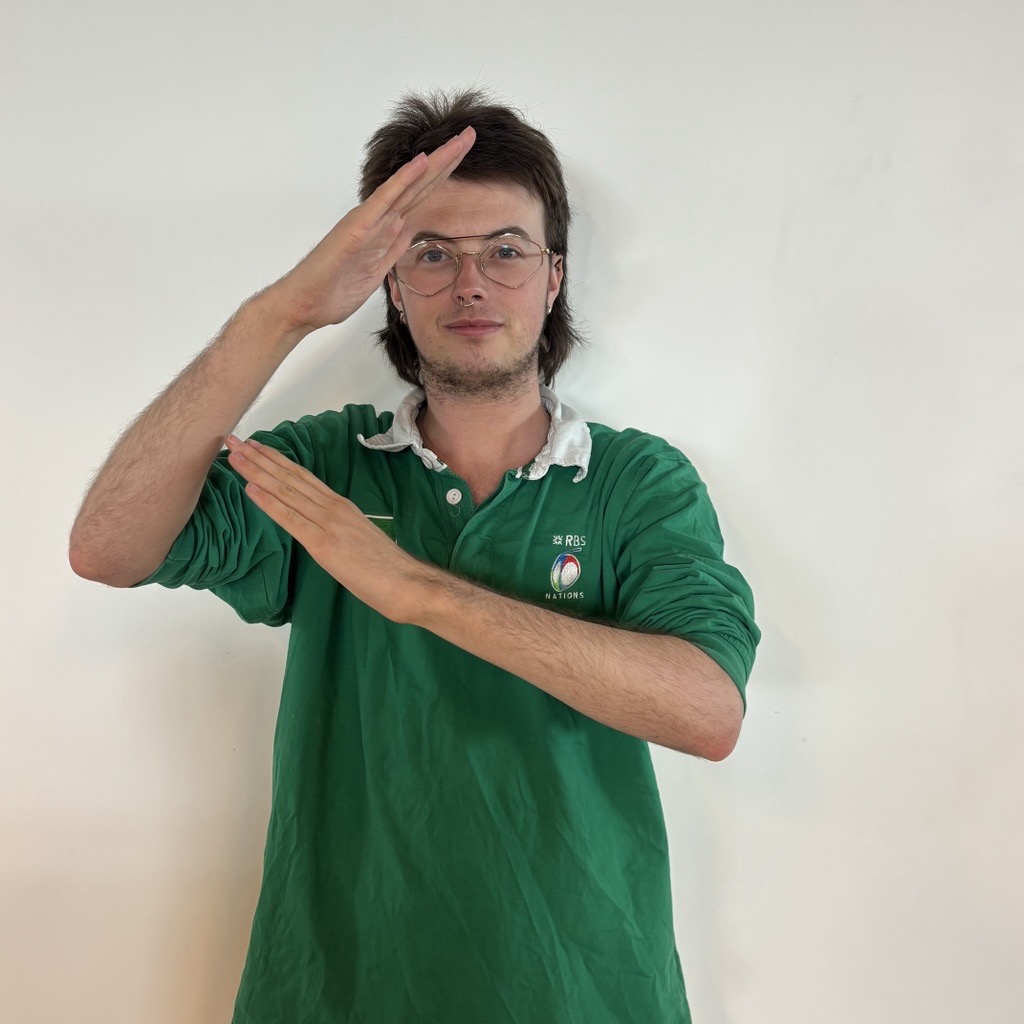}
\subcaption{Step 4.}
  \end{minipage}
  \caption{Square Routine Two}
  \label{sq2}
\end{figure}

Using both of these forms, we created a solid understanding of the shapes that could be made in this framework. The key difference between these was using the sides versus using the inside to map the shape. Both of these show right angles and perpendicularity but from alternate views- Often it is hard to see symmetries of a shape until you tilt your head and look at the object from a new perspective. The second clearly also includes a triangular form, showing that shapes can be built from others, not just from lines.

To merge all of these together and to create a unifying piece, we enacted a cannon of waves together, each student moving in rows. We conveyed to the students the importance of timing per the counts and each and every individual taking responsibility. Ultimately, this simple movement helped get across the concept of murmuration that we began with. We moved our arms in a bird-like way putting all of the above together in a fast rhythmic pattern, combining some movements from the squares.

Try this out yourself! 
 
\section{Post Dance}

After practising the routine, we asked for volunteers to show the rest of the room what they had just learnt. Indeed, we had a murmuration of enthusiastic children. One could clearly see the cannon effect we introduced, as well as the different shapes and symmetries that were embedded in the routine. We felt a sense of pride watching the pupils, particularly at their eagerness especially after having only around twenty minutes to practice. 

To bring more context to the movement,  we discussed the music with the composer of the piece, who talked through his experience of writing music, and how this is also interleaved with mathematics. He also played  live, which brought an extra dimension to the session that we were fortunate to have.

We wished to discuss the mathematics more concretely with the students, bringing in more formal notions of symmetry, but time was short. We have since realised that this is what is usually done in the classroom, and maybe a more efficient way would be to give the teachers a pack for a follow up session to what we conducted. 

We had a final surprise in store- we had recorded the performance of the students, and this was shown to the entire group to emphasise how local actions can greatly affect the global picture. Giving a sense of individual responsibility to the students, and showing them that yes, the big picture looks different to the individual picture but both are just as important (\'a la Gestalt). We hoped that despite this we conveyed to the students that perfection is not necessary for a great experience and therefore a great routine.

We made sure to leave space for a question and answer segment, giving the students an opportunity to lead, we were asked some interesting questions about our experiences and what led us to dance and mathematics. The best question was definitely the simple but effective `What is your favourite shape?', we can see that we made the mathematics of shape that little bit more personal!  

\section{Logistical comments}
A key back and forth we had was of the location of the session. From outside to inside to the shape of the lecture hall/designated area. Ultimately, we ended up in a curved lecture theatre, a semi-ideal space; this meant that the rows and columns of pupils was not always consistent. We anticipated to feel more of an affect from this, but the murmuration routine still happened fairly smoothly, which feels like a nod to the starlings only taking into account a subset of their peers whilst moving.

We also thought about doing this routine outside, as there was a more tiered set of steps, however, we ultimately decided that (being in England) the weather was too unpredictable and filtering almost 200 seven and eight year olds inside in case of rain may have been too disruptive. This decision was ultimately wise, as the enclosed space, at least for us, brought safety into enacting the movements freely without being seen by onlookers.

The dance movements were greatly inspired from the Paris 2024 Olympics and Paralympics handover dance routine shown at the closing ceremony of the Tokyo 2020 Paralympic Games. This was a routine of a dance form called tutting, which involves highly geometrical forms, choreographed by Sadeck Waff, and performed with 128 people, only half of these were professional dancers. This form of dancing is praised for its mesmerising qualities (which we find mathematics also holds) and how accessible it is in its nature. \cite{waffdance}

\section{Reflections and Conclusions}

Several days afterward, we received an email from a participating school, which stated that `I just wanted to thank you and all your colleagues for hosting the children. They came back absolutely full of it, so excited and keen to share what they had done. I have been learning the dance whilst sitting in my office!!' This shows to us that at least in the short term, the students found this experience memorable, and we hope that this has then embedded this topic in their minds. All in all, the  creative learning experience that we brought to these students was an extra layer to their education on these topics that they wouldn't usually have, particularly given how overstretched teachers are.

Another facilitator, upon watching and engaging in our session commented that he felt the sense of community behind the routine even just from watching it. What was experienced was `collective effervescence'- a sense of unity amongst the group (and spectators), which feels somewhat magical and akin to being swept away! \cite{Kronsted_2025} This coincides with the philosophy of Aristotle- the whole being more than the sum of its parts, but also with psychological Gestalt theories showing that  every part indeed counts to bring this feeling of collective unity.

Immediately after the session we wrote some comments on our own thoughts:
\begin{itemize}
    \item `Before this experience, I felt apprehensive about my own experience with dance/movement/sport, and how I often avoid such situations. Minutes into the first session, although I was stood in front of approximately 180 students, which is not something I do regularly, I began to feel a wash of calm. Being beside [the other facilitator] made the experience feel natural as we bounced ideas around the room, and I do wonder how different this session would look if only one of us (or someone else!) was doing it. 
    After the day was over, I felt motivated, the energy that was in that room was electric, although, I was exhausted! If we were to do this again, which I hope we can, I would attempt to change the order of the sessions. The students were very keen on the dancing part, which is exactly what we wanted, but it was then hard to get them to settle down afterward. This was particularly prominent in the afternoon session, as they had been working all morning on a different activity. It was therefore a lot more difficult than I expected to connect with the mathematics of what we were doing explicitly. The most obvious option to altering this is to show the prerequisites in mathematics before the dance routine. However, we run the risk of this becoming too much like a lecture.'

    \item `From my own perspective as a facilitator, the session began with a considerable sense of nervousness. Standing in front of such a large group of young people in a space that was not necessarily ideal for movement-based activities felt daunting, and I remember thinking, I have no idea what I am going to say. The scale of the audience, combined with uncertainties around choreography and logistics, initially created a sense of fear. However, once the session began, these concerns quickly diminished. The collaborative nature of the facilitation was crucial; by responding to one another in real time and building on each other's contributions, a natural rhythm emerged that made the experience feel both safe and connected. This sense of connection between facilitators appeared to influence the wider atmosphere of the room and helped establish rapport with the students. One particularly encouraging observation was the students' willingness to contribute during the opening discussion. Even when unsure of the \textbf{correct} answer, they eagerly shared ideas and drew on their existing knowledge of birds, movement, and the natural world. In doing so, the session benefited from a progression from known to unknown, grounding abstract mathematical ideas in familiar experiences. Seeing the ambassadors perform the routine was a pivotal moment; the visible excitement and curiosity in the students' faces demonstrated the power of modelling and aspiration. Although the choreography required significant spatial awareness and coordination, the students embraced the challenge enthusiastically. Watching the final performance brought an overwhelming sense of joy and pride, particularly given the relatively short amount of time available for practice. Reflecting afterwards, I felt that the mathematical ideas could have been made more explicit, either by extending the session or by introducing some of the mathematical concepts before the movement activities began. Nevertheless, the experience reinforced for me the value of embodied, collaborative learning, and highlighted how creativity, movement, and mathematics can come together to create memorable and meaningful educational experiences'.

\end{itemize}

The obvious next steps are to figure out how we can bring this to a classroom, or a non STEAM based event, as it is clear that the background of this event gave us a lot more flexibility than everyday teaching! We also wonder how to bring in the sense of randomness and the feeling of the routine being fleeting that starlings have during their murmurations; obviously they have more external factors to take into account, but how could we bring those to light in a different form? We think that a feeling of randomness could be brought into this routine, one way that we theorise is by letting a student choose the movement and this being propagated through the crowd- either actually randomly, or in a systematic way.

We wonder how this framework could benefit students who are disabled and/or neurodivergent, both inside and outside the classroom, as we recognise that this dance already has accessibility built into it.

Finally, a natural next step could be how to bring this to other mathematical concepts and also older students, as they may be more resistant to unconventional teaching styles.

\section{Acknowledgements}
The authors would like to thank Nick Elliott and the Widening Access Team at the University of Hertfordshire for the opportunity to present this session and for their help along the way. We would like to give thanks to all of the ambassadors who helped make this possible, to David Allen for his invaluable advice and feedback, and to Aryan Pise for the composition of the music used and exquisite violin playing for the students. Finally, we would like to thank the students and teachers from Yewtree Primary School, Fairlands Primary School, Howe Dell Primary School,
Oakview Primary School and 
St Vincent de Paul Catholic Primary School for their eager participation in the event. 

\bibliography{movement}
\bibliographystyle{plain}
\end{document}